\documentclass[a4paper,12pt]{article}
\usepackage{amsmath}
\usepackage{cite}

\begin{document}

\tolerance=500

\begin{center}
{\large\bf {Some extensions of Ramanujan's ${}_1\psi_1$ summation
formula}}
\end{center}

\begin{center}
    {N.M. Vildanov}\footnote{I.E. Tamm Department of Theoretical
    Physics, P.N. Lebedev Physics Institute, 119991 Moscow, Russia}
\end{center}

\begin{center}
\begin{minipage}{4.5in}
    \small {We have found several summation formulas that extend Ramanujan's psi
    sum. First contains a parameter $\alpha=1/N$, $N$ is a positive
    integer, and transforms to $q$-beta integral in the limit $N\to\infty$.
    The other is a $q$-analogue of generalized binomial theorem, expansion
    of $(1+t)^a$ in powers of $t^\alpha$, $0<\alpha\leq 1$, and it expresses
    the sum of a certain series in terms of a bibasic integral.}
\end{minipage}
\end{center}

\begin{flushleft}
    {\bf{1. Introduction}}
\end{flushleft}

Given a function $f_b(x)$ one can define the band limited function
$g_b(y)$ as follows
$$
    g_b(y)=\frac{1}{2\pi}\int_{-\pi}^{\pi}f_b(x)e^{ixy} dx.
$$
Then using Poisson summation formula, it can be shown that for
$0<\alpha\leq 1$
\begin{equation}\label{general_formula}
        \alpha\sum_{n=-\infty}^{\infty}g_b(y_1+\alpha n)g_c^*(y_2+\alpha n)=\frac{1}{2\pi}\int_{-\pi}^{\pi}f_b(x)f_c^*(x)e^{ix(y_1-y_2)}
        dx.
\end{equation}
Here $b$ denotes a set of parameters and, therefore, $g_b$ and
$g_c$ can be different functions. We note three special cases of
this identity:

(i) If $y_1=y_2=0$, $\alpha=1$, then we get Parseval's identity.

(ii) If $f_c(x)=e^{-ixy}$, then we get sampling theorem for a band
limited function $g(y)$ in the form

\begin{equation}\label{sampling}
    \sum_{n=-\infty}^{\infty}\frac{\sin\pi(y-\alpha n)}{\pi (y-\alpha n)}g(\alpha
    n)=\frac{1}{\alpha}~g(y)
\end{equation}

(iii) Substituting $f_b(x)=\cos^b(x/2)$ and using the value of the
integral\cite{gradshtein}
\begin{equation}\label{beta}
    \int_{0}^{\pi/2}\cos^b(x)\cos
ax=\frac{\pi\Gamma(b+1)}{2^{b+1}\Gamma\left(\frac{a+b}{2}+1\right)\Gamma\left(\frac{b-a}{2}+1\right)}
\end{equation}
we get generalization of the Dougall bilateral sum
\begin{align}\label{dougall}
    \nonumber &\alpha\sum_{n=-\infty}^{\infty}\frac{1}{\Gamma(a+\alpha n)\Gamma(b-\alpha n)\Gamma(c+\alpha n)\Gamma(d-\alpha
    n)}\\
    &\qquad {}
    =\frac{\Gamma(a+b+c+d-3)}{\Gamma(a+b-1)\Gamma(a+d-1)\Gamma(c+b-1)\Gamma(c+d-1)}.
\end{align}

There is a more general theorem, extension of Leibniz rule to
fractional derivatives\cite{osler1}. It states that if $D^a$
denotes a fractional derivative then
\begin{equation}\label{leibnize}
    D^a(uv)=\alpha\sum_{n=-\infty}^{\infty}\frac{\Gamma(a+1)}{\Gamma(a-c-\alpha n+1)\Gamma(\alpha
    n+c+1)}D^{c+\alpha n}uD^{a-c-\alpha n}v.
\end{equation}
Usefulness of this formula stems from the fact that many special
functions can be represented as fractional derivatives. In this
manner it is possible to compute sums of many complicated infinite
series \cite{osler1}.

Corresponding to formula \eqref{leibnize} is the following
generalization of binomial theorem\cite{osler2}
\begin{equation}\label{binomial_alpha}
    \alpha\sum_{n=-\infty}^{\infty}\frac{\Gamma(a+1)}{\Gamma(a-c-\alpha n+1)\Gamma(\alpha
    n+c+1)}t^{c+\alpha n}=(1+t)^a
\end{equation}
for $|t|=1$ ($t^{c+\alpha n}$ in \eqref{binomial_alpha} can be
replaced by $t^{a-c-\alpha n}$ by symmetry). The case $c=0$ is the
expansion of $(1+t)^a$ in powers of $t^\alpha$. Coefficients $A_n$
in the expansion of the function $f(x),~-\pi<x<\pi$ in powers of
$e^{i\alpha x}$,~$0<\alpha\leq 1$
$$
    f(x)=\sum_{n=-\infty}^{\infty}A_n e^{i\alpha xn}
$$
are given by
$A_n=\frac{\alpha}{2\pi}\int_{-\pi}^{\pi}f(x)e^{-i\alpha xn}dx$.
~One can use the integral \eqref{beta} to compute coefficients
$A_n$ for the function $(1+e^{ix})^a$ and see that the result
coincides with the case $c=0$ of \eqref{binomial_alpha}.

When $\alpha=1$, eq. \eqref{binomial_alpha} reduces
to\cite{riemann,heaviside,hardy,watanabe}
\begin{equation}\label{binomial_riemann}
    \sum_{n=-\infty}^{\infty}\frac{\Gamma(a+1)}{\Gamma(a-c-n+1)\Gamma(
    n+c+1)}t^{c+n}=(1+t)^a.
\end{equation}
This formula can be proved by writing
$(1+t)^a=(1+t)^{a-c}\cdot(1+t)^c$, expanding first factor in
powers of $t$ and the second factor in inverse powers, $1/t$, and
then computing the coefficient in front of $t^{c+n}$ by using
Gauss' sum for the hypergeometric function. Combining the formulas
for $0<\alpha\leq 1, c=0$ and $\alpha=1,c\neq 0$ one can obtain
the general formula \eqref{binomial_alpha}.

In the following, we will use the conventional notation for
$q$-products and basic hypergeometric series\cite{gasper_rahman}.
It is well known that Ramanujan's psi sum
\begin{align}\label{}
    \nonumber
    &{}_1\psi_1(a,b;q,z)=\sum_{n=-\infty}^\infty\frac{(a;q)_n}{(b;q)_n}z^n\\
    &=\frac{(q,b/a,az,q/az;q)_\infty}{(b,q/a,z,b/az;q)_\infty},\qquad
    |q|<1,~|b/a|<|z|<1
\end{align}
is a $q$-analogue of \eqref{binomial_riemann}, i.e. in the limit
$q\rightarrow 1^-$ it gives \eqref{binomial_riemann}, and for
$b=q$ it reduces to $q$-binomial theorem
\begin{equation}\label{}
    \nonumber
    \sum_{n=0}^\infty\frac{(a;q)_n}{(q;q)_n}z^n=\frac{(az;q)_\infty}{(z;q)_\infty},\qquad
    |z|<1.
\end{equation}

The aim of this paper is to obtain $q$-analogue of identity
\eqref{binomial_alpha}. It seems that the question of obtaining
$q$-analogues of sums like \eqref{binomial_alpha} or
\eqref{dougall} has not been considered in the literature. We note
also that basic hypergeometric series are not suitable to solve
this problem.

First, in sections 2 and 3, we will obtain some extensions of
Ramanujan's psi sum and Bailey's summation for the
very-well-poised ${}_6\psi_6$ containing an additional parameter
$\alpha =1/N$, $N$ is a positive integer. This will enable to
obtain an alternative derivation of some $q$-beta integrals as
discussed in section 4. Section 5 contains some new results on
summation of series which does not belong to the family of basic
hypergeometric series, eqs. \eqref{summation} and
\eqref{third_extension}. In the last section, we consider some
limiting cases of the identity \eqref{third_extension} and show
how it can be reduced to \eqref{binomial_alpha}.

\begin{flushleft}
    \large{\bf{2. First extension}}
\end{flushleft}

The bilateral series
\begin{equation}\label{first_sum_q}
    \sum_{n=-\infty}^{\infty}\frac{(bp^{n};q)_\infty}{(ap^{n};q)_\infty}z^n
\end{equation}
where $p=q^\alpha$, $|b/a|^\alpha<|z|<1$, is proportional to
${}_1\psi_1(a,b;q,z)$ when $\alpha=1$. In the limit $q\rightarrow
1^-$, it gives the infinite series
\begin{equation}\label{first_sum_gamma}
    \sum_{n=-\infty}^{\infty}\frac{\Gamma(a+\alpha n)}{\Gamma(b+\alpha
    n)}z^n.
\end{equation}
The series in eq.\eqref{binomial_alpha} can be, in principle,
obtained from \eqref{first_sum_gamma}. Now we will show how to
compute \eqref{first_sum_q} when $\alpha=1/N$, $N$ is a positive
integer.

We will utilize the method used for analytic continuation of basic
hypergeometric series\cite{watson,gasper_rahman}. Consider the
function
\begin{equation}\label{function1}
    \frac{(bq^{\alpha s};q)_\infty}{(aq^{\alpha
    s};q)_\infty}~\frac{\pi(-z)^s}{\sin\pi s}
\end{equation}
where $q=e^{-\omega}$, $p=e^{-\Omega}$, so $\Omega=\alpha\omega$.
We suppose that $0<q<1$, or equivalently $\omega, \Omega>0$. The
sum of the residues of \eqref{function1} at the points
$s=0,\pm1,\pm 2,...$ is equal to \eqref{first_sum_q} and at the
points $s=-{n}/{\alpha}+(\ln a+2\pi i m){\Omega^{-1}}$,
$n=0,1,2,3,...$, $m=0,\pm1,\pm 2,...$ to the double series
\begin{equation*}\label{}
    \sum_{n=0}^\infty\sum_{m=-\infty}^{\infty}\frac{(bq^{-n}/a;q)_\infty}{(q^{-n};q)_n (q;q)_\infty}~
    \frac{\pi(-z)^{-{n}/{\alpha}+(\ln a+2\pi i m){\Omega^{-1}}}}{\Omega\cdot\sin\pi[-{n}/{\alpha}+(\ln a+2\pi i
    m){\Omega^{-1}}]}.
\end{equation*}
If $1/\alpha$ is an integer, then the two sums can be separated.
The sum over $n$ can be computed using $q$-binomial theorem
\begin{align*}\label{}
    &\sum_{n=0}^\infty\frac{(bq^{-n}/a;q)_\infty}{(q^{-n};q)_n
    (q;q)_\infty}z^{-{n}N}\\
    &=\frac{(b/a;q)_\infty}{(q;q)_\infty}\sum_{n=0}^\infty\frac{(aq/b;q)_n}{(q;q)_n}(b/az^N)^n\\
    &=\frac{(b/a,q/z^N;q)_\infty}{(q,b/az^N;q)_\infty}.
\end{align*}
The sum over $m$ is\cite{gasper}
\begin{equation}\label{m_sum}
    \sum_{m=-\infty}^{\infty}\frac{\pi(-z)^{(\ln a+2\pi i m){\Omega^{-1}}}}{\Omega\cdot\sin(\pi\ln a/\Omega+2\pi^2 i
    m\Omega^{-1})}=-\frac{(p,p,az,p/az;p)_\infty}{(a,p/a,z,p/z;p)_\infty}.
\end{equation}
Since the sum of all of the residues of function \eqref{function1}
is zero we get the result
\begin{equation}\label{first_extension}
    \sum_{n=-\infty}^{\infty}\frac{(bp^{n};q)_\infty}{(ap^{n};q)_\infty}z^n=\frac{(b/a,q/z^N;q)_\infty}{(q,b/az^N;q)_\infty}\cdot
\frac{(p,p,az,p/az;p)_\infty}{(a,p/a,z,p/z;p)_\infty}, \qquad
p=q^{1/N}.
\end{equation}

Before proceeding further we note a connection of the formula
\eqref{first_extension} with theta functions. This connection is
not surprizing in view of the fact that Jacobi triple product is a
special case of Ramanujan's sum.

Replacing $n$ in eq. \eqref{first_extension} by $k+nN$, where
$k=0,1,2,...,(N-1)$, $~n=0,\pm 1,\pm 2,...$ we get
\begin{align*}\label{}
    &\sum_{n=-\infty}^\infty\sum_{k=0}^{N-1}\frac{(bq^{n+k/N};q)_\infty}{(aq^{n+k/N};q)_\infty}z^{k+nN}\\
    &=\sum_{k=0}^{N-1}z^k\frac{(bq^{k/N};q)_\infty}{(aq^{k/N};q)_\infty}\sum_{n=-\infty}^\infty \frac{(aq^{k/N};q)_n}{(bq^{k/N};q)_n}
    z^{k+nN}\\
    &=\sum_{k=0}^{N-1}z^k\frac{(bq^{k/N};q)_\infty}{(aq^{k/N};q)_\infty}\cdot \frac{(q,b/a,aq^{k/N}z^N,q^{1-k/N}/az^N;q)_\infty}
    {(bq^{k/N},q^{1-k/N}/a,z^N,b/az^N;q)_\infty}.
\end{align*}
So, we obtain a relation between infinite products
\begin{align*}\label{}
    \sum_{k=0}^{N-1}&z^k\frac{(aq^{k/N}z^N,q^{1-k/N}/az^N;q)_\infty}{(aq^{k/N},q^{1-k/N}/a;q)_\infty}\\
    &=\frac{(z^N,q/z^N;q)_\infty}{(q,q;q)_\infty}\cdot
    \frac{(p,p,az,p/az;p)_\infty}{(a,p/a,z,p/z;p)_\infty},\qquad
    p=q^{1/N}.
\end{align*}
When expressed in terms of theta functions this identity reads
\begin{align*}
    &\sum_{k=0}^{N-1}e^{2iky}\frac{\theta_1(x+Ny+\pi ik \tau|N\tau)}{\theta_1(x+\pi
    ik\tau|N\tau)}=\frac{\theta_1^\prime(0|\tau)}{\theta_1^\prime(0|N\tau)}\frac{\theta_1(Ny|N\tau)\theta_1(x+y|\tau)}{\theta_1(x|\tau)\theta_1(y|\tau)}.
\end{align*}

When $N=2$ eq. \eqref{first_extension} can be reduced to
\begin{align*}
            (aq^{\frac{1}{2}},{q^{\frac{1}{2}}}/{a},az^2,q/az^2;q)_\infty+z
            (a,q/{a},az^2q^{\frac{1}{2}},q^{\frac{1}{2}}/az^2;q)_\infty=\\
        (q^{\frac{1}{2}};q)^2_\infty(-z,-q^{\frac{1}{2}}/z,az,q^{\frac{1}{2}}/az;q^{\frac{1}{2}})_\infty.
\end{align*}
This identity can be written in terms of theta functions, after
defining new variables, as follows
\begin{equation*}\label{}
    \theta_3(x|\tau)\theta_4(y|\tau)=\theta_4(x+y|2\tau)\theta_4(y-x|2\tau)+\theta_1(x+y|2\tau)\theta_1(y-x|2\tau).
\end{equation*}
When $x=y$ this is Landen's transform\cite{whittaker}
\begin{equation*}\label{}
    \frac{\theta_3(x|\tau)\theta_4(x|\tau)}{\theta_4(2x|2\tau)}=\frac{\theta_3(0|\tau)\theta_4(0|\tau)}{\theta_4(0|2\tau)}.
\end{equation*}

\begin{flushleft}
    \large{\bf{3. Second extension}}
\end{flushleft}

The following consequence of Ramanujan's psi sum
\begin{equation*}\label{}
    \sum_{n=-\infty}^{\infty}\frac{(-x)^n}{1-q^{\beta+n}}=
    \frac{(q,q;q)_\infty}{(q^{\beta},q^{1-\beta};q)_\infty}\frac{(-q^\beta x,-q^{1-\beta}/x;q)_\infty}{(-x,-q/x;q)_\infty}
\end{equation*}
is the $q$-analogue of the identity
$$
    \frac{\pi x^\beta}{\sin\pi \beta}=\sum_{-\infty}^\infty \frac{(-x)^n}{\beta
    -n}.
$$
In view of this observation it is more proper to consider the sum
\begin{equation}\label{second_sum}
    \sum_{n=-\infty}^\infty\frac{(bq^{\alpha n},q^{1-\alpha n}/a;q)_\infty}{(-axq^{\alpha n},-q^{1-\alpha n}/ax;q)_\infty}
\end{equation}
instead of \eqref{first_sum_q}. This sum is proportional to
${}_1\psi_1(a,b;q,z)$ when $\alpha=1$ and it converges in the same
ring $|b/a|<|z|<1$ as ${}_1\psi_1(a,b;q,z)$ does for arbitrary
$0<\alpha\leq 1$.

As in the previous section, we will compute the sum for
$\alpha=1/N$, where $N$ is a positive integer. Consider the
function
\begin{equation}\label{function2}
    \frac{(bq^{\alpha s},q^{1-\alpha s}/a;q)_\infty}{(-axq^{\alpha s},-q^{1-\alpha s}/ax;q)_\infty}~\frac{\pi(-y)^s}{\sin\pi s}
\end{equation}
where a new variable has been introduced, $y$, which is absent in
\eqref{second_sum}. This is an auxiliary variable and the limit
$y\rightarrow 1$ will be taken at the end of calculations.

Sum of the residues of function \eqref{function2} at $s=0,\pm
1,\pm 2,\pm 3,...$ is
\begin{equation*}\label{}
    \sum_{n=-\infty}^\infty\frac{(bq^{\alpha n},q^{1-\alpha n}/a;q)_\infty}{(-axq^{\alpha n},-q^{1-\alpha
    n}/ax;q)_\infty}~y^n
\end{equation*}
and at $s={n}/{\alpha}+[\ln (-ax)+2\pi i m]{\Omega^{-1}}$,
$n,m=0,\pm 1,\pm 2,...$ is
\begin{equation*}\label{}
    \sum_{n=0}^\infty\sum_{m=-\infty}^{\infty}\frac{(-bq^{-n}/ax,-xq^{1+n};q)_\infty}{(q^{-n},q^{1+n};q)_n (q;q)_\infty}~
    \frac{\pi(-y)^{-{n}/{\alpha}+[\ln (-ax)+2\pi i m]{\Omega^{-1}}}}{\Omega\sin\pi\{-{n}/{\alpha}+[\ln (-ax)+2\pi i
    m]{\Omega^{-1}}\}}
\end{equation*}
\begin{equation*}\label{}
    -\sum_{n=1}^\infty\sum_{m=-\infty}^{\infty}\frac{(-bq^{n}/ax,-xq^{1-n};q)_\infty}{(q^{n},q^{1-n};q)_{n-1} (q;q)_\infty}~
    \frac{\pi(-y)^{{n}/{\alpha}+[\ln (-ax)+2\pi i m]{\Omega^{-1}}}}{\Omega\sin\pi\{{n}/{\alpha}+[\ln (-ax)+2\pi
    im]{\Omega^{-1}}\}}
\end{equation*}
\begin{align*}\label{}
    &=\frac{(-b/ax,-xq;q)_\infty}{(q,q;q)_\infty}\sum_{n=-\infty}^\infty\frac{(-1/x;q)_n}{(-b/ax;q)_n}(-xy^N)^n\\
    &~~~\cdot\sum_{m=-\infty}^{\infty}\frac{\pi(-y)^{[\ln (-ax)+2\pi im]{\Omega^{-1}}}}{\Omega\sin[\pi\ln
    (-ax)/\Omega+2\pi^2
    im\Omega^{-1}]}\\
    &=-\frac{(b/a,y^N,q/y^N;q)_\infty}{(q,-xy^N,-b/axy^N;q)_\infty}\frac{(p,p,-axy,-p/axy;p)_\infty}{(-ax,-p/ax,y,p/y;p)_\infty}.
\end{align*}
Equating the sum of the residues to zero gives the identity
\begin{align}\label{y_identity}
    &\nonumber\sum_{n=-\infty}^\infty\frac{(bq^{\alpha n},q^{1-\alpha n}/a;q)_\infty}{(-axq^{\alpha n},-q^{1-\alpha
    n}/ax;q)_\infty}~y^n\\
    &=\frac{(b/a,y^N,q/y^N;q)_\infty}{(q,-xy^N,-b/axy^N;q)_\infty}\frac{(p,p,-axy,-p/axy;p)_\infty}{(-ax,-p/ax,y,p/y;p)_\infty},\quad
    p=q^{1/N}.
\end{align}
Taking the limit $y\to 1$ one obtains
\begin{equation*}
    \alpha \sum_{n=-\infty}^\infty\frac{(bq^{\alpha n},q^{1-\alpha n}/a;q)_\infty}{(-axq^{\alpha n},-q^{1-\alpha
    n}/ax;q)_\infty}=\frac{(q,b/a;q)_\infty}{(-x,-b/ax;q)_\infty},\qquad
    \alpha=1/N,
\end{equation*}
and redefining the parameters the result
\begin{equation}\label{alpha_psi_sum}
    \alpha \sum_{n=-\infty}^\infty\frac{(bq^{\alpha n},q^{1-\alpha n}/a;q)_\infty}{(xq^{\alpha n},q^{1-\alpha
    n}/x;q)_\infty}=\frac{(q,b/a;q)_\infty}{(x/a,b/x;q)_\infty}
\end{equation}
where $|b|<|x|<|a|$,$~\alpha=1/N$.

One may hope to extend this expression to $0<\alpha\leq 1$ by
Carlson's theorem. However, this is not possible.

\begin{flushleft}
    \large{\bf{4. Application to $q$-beta integrals}}
\end{flushleft}

There are some consequences of the formulas obtained in the
previous section. In the limit $\alpha\to 0$, summation over $n$
in \eqref{alpha_psi_sum} can be replaced by integration and we
obtain the $q$-beta integral\cite{askey,gasper}
\begin{equation*}\label{beta_integral}
    \int_0^\infty\frac{(bt,q/at;q)_\infty}{(-t,-q/t;q)_\infty}\frac{dt}{t}=\frac{(q,b/a;q)_\infty}{(-1/a,-b;q)_\infty}\ln\frac{1}{q}.
\end{equation*}
Since $N$ does not enter the right hand side of
\eqref{alpha_psi_sum}, one can see how Ramanujan's psi sum
transforms to $q$-beta integral when $N$ increases from 1 to
$\infty$.

Analogously, the limit $N\to\infty$ of eq.\eqref{y_identity} is
the integral\cite{ramanujan,askey,gasper}
\begin{equation*}\label{beta_integral2}
    \int_0^\infty\frac{(-q^bt,-q^{1-a}/t;q)_\infty}{(-t,-q/t;q)_\infty}t^{c-1}dt=
    \frac{\pi}{\sin\pi
    c}\frac{(q^{b-a},q^c,q^{1-c};q)_\infty}{(q,q^{b-c},q^{c-a};q)_\infty}.
\end{equation*}

Eq. \eqref{alpha_psi_sum} can be derived without using countour
integration. We will demonstrate the method showing how one can
obtain Askey's $q$-beta integral \cite{askey2} by generalizing
Bailey's summation for ${}_6\psi_6$. The reason why this method
works is the same as in other derivations of $q$-beta integrals
(see, e.g., \cite{ismail,suslov}).

Bailey's formula is \cite{gasper_rahman}
\begin{align}\label{bailey}
    \nonumber &{}_6\psi_6 \left[\begin{array}{c}
      qa^{\frac{1}{2}}, -qa^{\frac{1}{2}},b,c,d,e \\
      a^{\frac{1}{2}}, -a^{\frac{1}{2}},aq/b,aq/c,aq/d,aq/e \\
    \end{array};q,\frac{qa^2}{bcde}\right]\\
    &=\frac{(aq,aq/bc,aq/bd,aq/be,aq/cd,aq/ce,aq/de,q,a/q;q)_\infty}{(aq/b,aq/c,aq/d,aq/e,q/b,q/c,q/d,q/e,qa^2/bcde;q)_\infty}.
\end{align}
One can compute more general sum
\begin{align}\label{sn}
    \nonumber S_N=\alpha \sum_{n=-\infty}^\infty \frac{(abq^{\alpha n},acq^{\alpha n},adq^{\alpha n},
    aeq^{\alpha n};q)_\infty}{(
    aq^{1+\alpha n},-aq^{1+\alpha n},q^{1-\alpha n}/a,-q^{1-\alpha n}/a;q)_\infty}\\
    \cdot \frac{(bq^{-\alpha n}/a,cq^{-\alpha n}/a,dq^{-\alpha n}/a,eq^{-\alpha n}/a;q)_\infty}
    {(aq^{\frac{1}{2}+\alpha n},-aq^{\frac{1}{2}+\alpha n},q^{\frac{1}{2}-\alpha n}/a,
    -q^{\frac{1}{2}-\alpha n}/a;q)_\infty}
\end{align}
where $\alpha=1/N$. $S_1$ is
\begin{equation*}
    S_1=\frac{(bc/q,bd/q,be/q,cd/q,ce/q,de/q,q;q)_\infty}{(bcde/q^3;q)_\infty}
\end{equation*}
by \eqref{bailey} and does not depend on $a$. We can sum in
\eqref{sn} over $k+nN$, $k=0,1,2,...,N-1$ instead of $n$. This
leads to $N$ sums like $S_1$ with $a$ replaced by $aq^{k/N}$.
Since $S_1$ does not depend on $a$, summation over $k$ gives
$\frac{1}{N}\sum_{k=0}^{N-1}S_1=S_1$, and we obtain
\begin{equation*}
    S_N=S_1,\quad N=1,2,3,...
\end{equation*}
Substituting $a=i$ and taking the limit $N\to \infty$ one obtains
the integral
\begin{equation*}
    S_\infty=\int_{0}^\infty \frac{(ibu,-ib/u,icu,-icu,idu,id/u,ieu,-ie/u;q)_\infty}{(iq^{\frac{1}{2}}u,-iq^{\frac{1}{2}}/u,
    iqu,-iq/u,-iq^{\frac{1}{2}}u,iq^{\frac{1}{2}}/u,
    -iqu,iq/u;q)_\infty}\frac{du}{u\ln q^{-1}}.
\end{equation*}
So, we have proved Askey's theorem which in our notation can be
stated as $S_\infty=S_1$.

\begin{flushleft}
    \large{\bf{5. Third extension}}
\end{flushleft}

We can also consider the function defined by infinite series
\begin{equation}\label{f}
    f(a,b,z)=\sum_{n=-\infty}^\infty (bq^n,q^{1-n}/a;p)_\infty
    (-z)^n q^{n(n-1)/2}
\end{equation}
where $q>p$, so one can write $q=p^\alpha, 0<\alpha<1$; or $q=p$
and $|b|<|z|<|a|$. Here, we have reversed the notation to draw
some parallels with the work \cite{stanton} where the following
bibasic sum has been computed
\begin{align}\label{bibasic_sum}
    \nonumber
    &\frac{(q,qa^2;q)_\infty}{(qae^{i\theta},qae^{-i\theta};q)_\infty}(be^{i\theta},be^{-i\theta};p)_\infty\\
    &=\sum_{k=0}^\infty\frac{1-a^2q^{2k}}{1-a^2}\frac{(a^2,ae^{i\theta},ae^{-i\theta};q)_k}{(q,qae^{i\theta},qae^{-i\theta};q)_k}
    (-1)^kq^{k(k+1)/2}(abq^k,bq^{-k}/a;p)_\infty
\end{align}
with the same restriction on parameters $p,q,a,b$ as above.

There is a simple transformation formula
\begin{equation*}
    f(a,b,z)=-zf(1/b,1/a,1/x)
\end{equation*}
which can be extracted directly from the definition \eqref{f}. By
considering the poles of the function
\begin{equation*}
    \frac{(bq^{s},q^{1-s}/a;p)_\infty}{(xp^{s},p^{1-s}/x;p)_\infty}~\frac{\pi(-z)^s}{\sin\pi
    s},
\end{equation*}
one can derive the transformation formula
\begin{equation}\label{transformation}
    f(a,b,z)=\frac{(z,q/z;q)_\infty}{(z/y,qy/z;q)_\infty}f(a/y,b/y,z/y)
\end{equation}
which means that $f(a,b,z)/(z,q/z;q)_\infty$ depends in fact only
on two variables, e.g. $z/a$ and $b/z$. We will show this also by
finding a closed expression for $f(a,b,z)$.

Transformation formula \eqref{transformation} shows that $z=q^m$,
$m$ is an integer, are zeroes of $f(a,b,z)$. It can be seen also
from the definition of $f$ that $f(c,1/c,1)=0$.

Substituting $y=\sqrt{ab}$ into \eqref{transformation} and using
the fact $f(c,1/c,1)=0$, one can derive the identity
\begin{align*}
    &f(a,b,\sqrt{ab})\\
    &=\frac{(\sqrt{ab},q/\sqrt{ab};q)_\infty}{(q,q;q)_\infty}\sum_{n=1}^\infty
    (\sqrt{b/a}q^n,\sqrt{b/a}q^{1-n};p)_\infty (-1)^{n-1}(2n-1)
    q^{n(n-1)/2}.
\end{align*}
which is generalization of the identity
\begin{equation*}
    (q;q)^3_\infty=\sum_{n=1}^\infty
(-1)^{n-1}(2n-1)
    q^{n(n-1)/2}
\end{equation*}
known from the theory of theta functions.

Now we will show how to find a closed expression for
$f(aq/p,b,z)$. First, we compute the sum
\begin{align*}
    \nonumber & f(aq/p,a,z)=\sum_{n=-\infty}^\infty (aq^n,pq^{-n}/a;p)_\infty (-z)^n
    q^{n(n-1)/2}.
\end{align*}
We use Jacobi's triple product identity to expand the infinite
product $~(aq^n,pq^{-n}/a;p)_\infty$ and then to resum reversing
the order of summation
\begin{align*}
    \nonumber & f(aq/p,a,z)=\frac{1}{(p;p)_\infty}\sum_{n=-\infty}^\infty \sum_{m=-\infty}^\infty
    (-aq^n)^m
    p^{m(m-1)/2} (-z)^n
    q^{n(n-1)/2}\\
    &=\frac{(q;q)_\infty}{(p;p)_\infty}\sum_{m=-\infty}^\infty
    (-a)^m
    p^{m(m-1)/2} (zq^m,q^{1-m}/z;q)_\infty\\
    &=\frac{(q;q)_\infty}{(p;p)_\infty}\sum_{m=-\infty}^\infty
    (a/z)^m
    (p/q)^{m(m-1)/2} (z,q/z;q)_\infty\\
    &= \frac{(q;q)_\infty(p/q;p/q)_\infty}{(p;p)_\infty}(-a/z,-pz/aq;p/q)_\infty
    (z,q/z;q)_\infty.
\end{align*}
Thus, we have established the summation formula
\begin{align}\label{summation}
    \nonumber &\sum_{n=-\infty}^\infty (aq^n,pq^{-n}/a;p)_\infty (-z)^n q^{n(n-1)/2}\\
    &=\frac{(q;q)_\infty(p/q;p/q)_\infty}{(p;p)_\infty}(-a/z,-pz/aq;p/q)_\infty
    (z,q/z;q)_\infty,\qquad q>p.
\end{align}

Now we will use \eqref{summation} to compute $f(aq/p,b,z)$. One
can use Ramanujan's psi sum to show that
\begin{equation*}
    (bq^n,pq^{-n}/a;p)_\infty=(p,b/a;p)_\infty\int_{-\pi}^\pi \frac{(yq^ne^{i\theta},pq^{-n}e^{-i\theta}/y;p)_\infty}
    {(ye^{i\theta}/a,be^{-i\theta}/y;p)_\infty}\frac{d\theta}{2\pi}.
\end{equation*}
So, one can write
\begin{align*}
    &f(aq/p,b,z)=\sum_{n=-\infty}^\infty (bq^n,pq^{-n}/a;p)_\infty (-z)^n
    q^{n(n-1)/2}\\
    &=\sum_{n=-\infty}^\infty (p,b/a;p)_\infty\int_{-\pi}^\pi \frac{(yq^ne^{i\theta},pq^{-n}e^{-i\theta}/y;p)_\infty}
    {(ye^{i\theta}/a,be^{-i\theta}/y;p)_\infty}\frac{d\theta}{2\pi} (-z)^n
    q^{n(n-1)/2}\\
    &=(p,b/a;p)_\infty
    \int_{-\pi}^\pi\frac{f(qye^{i\theta}/p,ye^{i\theta},z)}{(ye^{i\theta}/a,be^{-i\theta}/y;p)_\infty}\frac{d\theta}{2\pi},
\end{align*}
and finally
\begin{align}\label{third_extension}
    \nonumber & \sum_{n=-\infty}^\infty (bq^n,pq^{-n}/a;p)_\infty (-z)^n
    q^{n(n-1)/2}\\
    &=(q;q)_\infty (p/q;p/q)_\infty (b/a;p)_\infty (z,q/z;q)_\infty
    \int_{-\pi}^\pi
    \frac{(ye^{i\theta}/z,pze^{-i\theta}/qy;p/q)_\infty}{(-ye^{i\theta}/a,-be^{-i\theta}/y;p)_\infty}\frac{d\theta}{2\pi}.
\end{align}

Thus, we have expressed $f(aq/p,b,z)$ as a bibasic integral.
Interestingly, a bibasic extension of Nasrallah-Rahman integral
has been found recently in \cite{ismrah} with the use of the sum
\eqref{bibasic_sum}. This fact was helpful in establishing the
result \eqref{third_extension}, though the integral
\eqref{third_extension} is quite different from the integrals
considered in \cite{ismrah}.

It is easy to see that transformation formula
\eqref{transformation} follows from \eqref{third_extension} by
assigning different values to $y$.

\begin{flushleft}
    \large{\bf{6.}}
\end{flushleft}

Let us consider some limiting cases of the identity
\eqref{third_extension}. The limit $p\to q^-$ can be found by
noticing that according to Jacobi triple product identity and
Poisson summation formula one has
\begin{align*}
    &(q,-e^{i\theta},-qe^{-i\theta};q)_\infty\\
    &=\sum_{m=-\infty}^\infty
    e^{im\theta}q^{m(m-1)/2}\to \sum_{k=-\infty}^\infty
    \delta(\theta/2\pi-k), \quad q\to 1^-
\end{align*}
where $\delta$ is delta-function. So the integral can be trivially
computed and we obtain Ramanujan's psi sum.

The limit
\begin{equation}\label{limit}
    q=p^\alpha, ~0<\alpha<1, ~|z|=1, ~p\to 1^-
\end{equation}
is quite nontrivial. As a first step we rewrite
\eqref{third_extension} as
\begin{align}\label{third_extension2}
    \nonumber & \sum_{n=-\infty}^\infty \frac{(p^{b+\alpha n},p^{1-a-\alpha n};p)_\infty}{(p,p^{b-a};p)_\infty} e^{it\alpha
    n}
    q^{n(n-1)/2}\\
    &=\frac{(q;q)_\infty (p/q;p/q)_\infty}{(p;p)_\infty}
    (-e^{it\alpha},-qe^{-it\alpha};q)_\infty
    \int_{-\pi}^\pi \frac{(-e^{i\theta-it\alpha},-pe^{-i\theta+it\alpha}/q;p/q)_\infty}{(-e^{i\theta}p^{-a},-p^be^{-i\theta};p)_\infty}\frac{d\theta}{2\pi}
\end{align}
where $t$ is real. In the limit \eqref{limit}, left hand side of
this equation becomes
\begin{equation*}
    \sum_{n=-\infty}^{\infty}\frac{\Gamma(b-a)}{\Gamma(b+\alpha n)\Gamma(1-a-\alpha
    n)}e^{it\alpha n}.
\end{equation*}
Then we present the integrand in the form
\begin{equation*}
    \frac{(-e^{i\theta-it\alpha},-pe^{-i\theta+it\alpha}/q;p/q)_\infty}{(-e^{i\theta},-pe^{-i\theta};p)_\infty}\cdot
    \frac{(-e^{i\theta},-pe^{-i\theta};p)_\infty}{(-e^{i\theta}p^{-a},-p^be^{-i\theta};p)_\infty}.
\end{equation*}
We know that\cite{gasper_rahman} the limit
\begin{align*}
    &\lim_{p\to1^-}\frac{(-e^{i\theta},-pe^{-i\theta};p)_\infty}{(-e^{i\theta}p^{-a},-p^be^{-i\theta};p)_\infty}\\
    &=(1+e^{i\theta})^{-a}(1+e^{-i\theta})^{b-1}=(1+e^{i\theta})^{b-a-1}e^{i\theta (1-b)}
\end{align*}
exists and it is a smooth function when $b\geq a+1$. Now we can
rewrite the limit \eqref{limit} of the right hand side of
\eqref{third_extension2} in the symmetrical form
\begin{align*}
    \int_{-\pi}^\pi g(\theta,t) (1+e^{i\theta})^{b-a-1}e^{i\theta
    (1-b)}\frac{d\theta}{2\pi}
\end{align*}
where
\begin{align*}
     g(\theta,t)
    =\lim_{p\to 1^-}\frac{(p^\alpha,-e^{it\alpha},-p^\alpha e^{-it\alpha};p^\alpha)_\infty
    (p^{1-\alpha},-e^{i\theta-it\alpha},-p^{1-\alpha}e^{-i\theta+it\alpha};p^{1-\alpha})_\infty}
    {(p,-e^{i\theta},-pe^{-i\theta};p)_\infty}.
\end{align*}
The function $g(\theta,t)$ does not depend on $a$ and $b$. Thus,
the limit \eqref{limit} of eq. \eqref{third_extension2} is
\begin{equation}\label{limiting}
    \sum_{n=-\infty}^{\infty}\frac{\Gamma(b-a)}{\Gamma(b+\alpha n)\Gamma(1-a-\alpha
    n)}e^{it\alpha n}=\int_{-\pi}^\pi g(\theta,t) (1+e^{i\theta})^{b-a-1}e^{i\theta
    (1-b)}\frac{d\theta}{2\pi}.
\end{equation}
Substituting $b=a+1$ we get
\begin{equation*}
    \sum_{n=-\infty}^{\infty}\frac{\sin\pi(a+\alpha n)}{\pi(a+\alpha n)}e^{it\alpha n}=\int_{-\pi}^\pi g(\theta,t) e^{-ia\theta
    }\frac{d\theta}{2\pi}.
\end{equation*}
One can use \eqref{sampling} or calculate the sum directly by
other means to show that this is equivalent to
$$
    \int_{-\pi}^\pi g(\theta,t) e^{-ia\theta
    }\frac{d\theta}{2\pi}=\frac{1}{\alpha}e^{-iat},
$$
i.e. $\alpha g(\theta,t)$ acts like delta-function. Substituting
this back into \eqref{limiting} we obtain the binomial theorem
\eqref{binomial_alpha}. It is not clear whether consideration of
the limit \eqref{limit} of the identity \eqref{third_extension}
can be carried out more consistently and rigorously.

\begin{flushleft}
    \large{\bf{7. Concluding remarks}}
\end{flushleft}

One may notice that the series \eqref{alpha_psi_sum} have some
similarity with elliptic hypergeometric series\cite{spiridonov}.
However, the infinite product $(q^n;p)_\infty$ can not be obtained
as a finite combination of theta functions. Moreover, there is no
relation between the bases $p$ and $q$, like $q>p$, in the theory
of elliptic hypergeometric functions. It is not clear what is the
precise relation of the formulas presented in this paper to the
theory of elliptic hypergeometric series.

\end{document}